\documentclass[runningheads]{llncs}

\usepackage{fancybox}
\usepackage{pstricks}
\usepackage{amsfonts}
\usepackage{graphicx}
\usepackage{algorithm}
\usepackage{algorithmic}
\usepackage{amsmath}
\usepackage{amssymb}
\usepackage[square,comma,sort]{natbib}

\titlerunning{An Approach for the QR Factorization on Multicore Architectures}

\pdfoutput=1
\relax
\title{Parallel Tiled QR Factorization for Multicore Architectures}

\author{Alfredo Buttari\inst{1},
  Julien Langou\inst{3},
  Jakub Kurzak\inst{1},
  Jack Dongarra\inst{1}\inst{2}
  }

\institute{Department of Electrical Engineering and Computer
Science, University Tennessee, Knoxville, Tennessee \and Computer
Science and Mathematics Division, Oak Ridge National Laboratory, Oak
Ridge, Tennessee \and Department of Mathematical Sciences,
University of Colorado at Denver and Health Sciences Center,
Colorado}

\begin{document}
\maketitle
\setcounter{page}{1}

\begin{abstract}
As multicore systems continue to gain ground in the High Performance
Computing world, linear algebra algorithms have to be reformulated
or new algorithms have to be developed in order to take advantage of
the architectural features on these new processors. Fine grain
parallelism becomes a major requirement and introduces the necessity
of loose synchronization in the parallel execution of an operation.
This paper presents an algorithm for the QR factorization where the
operations can be represented as a sequence of small tasks that
operate on square blocks of data. These tasks can be dynamically
scheduled for execution based on the dependencies among them and on
the availability of computational resources. This may result in an
out of order execution of the tasks which will completely hide the
presence of intrinsically sequential tasks in the factorization.
Performance comparisons are presented with the LAPACK algorithm for
QR factorization where parallelism can only be exploited at the
level of the BLAS operations.
\end{abstract}

\section{Introduction}

In the last twenty years, microprocessor manufacturers have been
driven towards higher performance rates only by the exploitation of
higher degrees of {\em Instruction Level Parallelism} (ILP). Based
on this approach, several generations of processors have been built
where clock frequencies were higher and higher and pipelines were
deeper and deeper. As a result, applications could benefit from
these innovations and achieve higher performance simply by relying
on compilers that could efficiently exploit ILP. Due to a number of
physical limitations (mostly power consumption and heat dissipation)
this approach cannot be pushed any further. For this reason, chip
designers have moved their focus from ILP to {\em Thread Level
  Parallelism} (TLP) where higher performance can be achieved by
replicating execution units (or {\em cores}) on the die while keeping
the clock rates in a range where power consumption and heat
dissipation do not represent a problem. Multicore processors clearly
represent the future of computing. It is easy to imagine that
multicore technologies will have a deep impact on the High Performance
Computing (HPC) world where high processor counts are involved
and, thus,  limiting power consumption and heat dissipation is a major
requirement. The Top500~\cite{top500} list released in June 2007 shows that
the number of dual-core Intel Woodcrest processors grew in six months
(i.e. from the previous list) from 31 to 205 and that 90 more systems
are based on dual-core AMD Opteron processors.

Even if many attempts have been made in the past to develop
parallelizing compilers, they proved themselves efficient only on a
restricted class of problems. As a result, at this stage of the
multicore era, programmers cannot rely on compilers to take advantage
of the multiple execution units present on a processor. All the
applications that were not explicitly coded to be run on parallel
architectures must be rewritten with parallelism in mind. Also,
those applications that could exploit parallelism may need considerable
rework in order to take advantage of the fine-grain parallelism
features provided by multicores.

The current set of multicore chips from Intel and AMD are for the most
part multiple processors glued together on the same chip. There are
many scalability issues to this approach and it is unlikely that type
of architecture will scale up beyond 8 or 16 cores. Even though it is
not yet clear how chip designers are going to address these issues, it
is possible to identify some properties that algorithms must have in
order to match high degrees of TLP:

\begin{description}
\item[fine granularity:] cores are (and probably will be) associated
  with relatively small local memories (either caches or explicitly managed
  memories like in the case of the STI Cell~\cite{isscc_2005_cell_desing} architecture or the Intel
  Polaris\cite{polaris} prototype). This requires splitting an operation into tasks
  that operate on small portions of data in order to reduce bus
  traffic and improve data locality.
\item[asynchronicity:] as the degree of TLP grows and granularity of
  the operations becomes smaller, the presence of synchronization
  points in a parallel execution seriously affects the efficiency of
  an algorithm.
\end{description}

\begin{figure}[!h]
  \begin{center}
    \includegraphics[width=\textwidth]{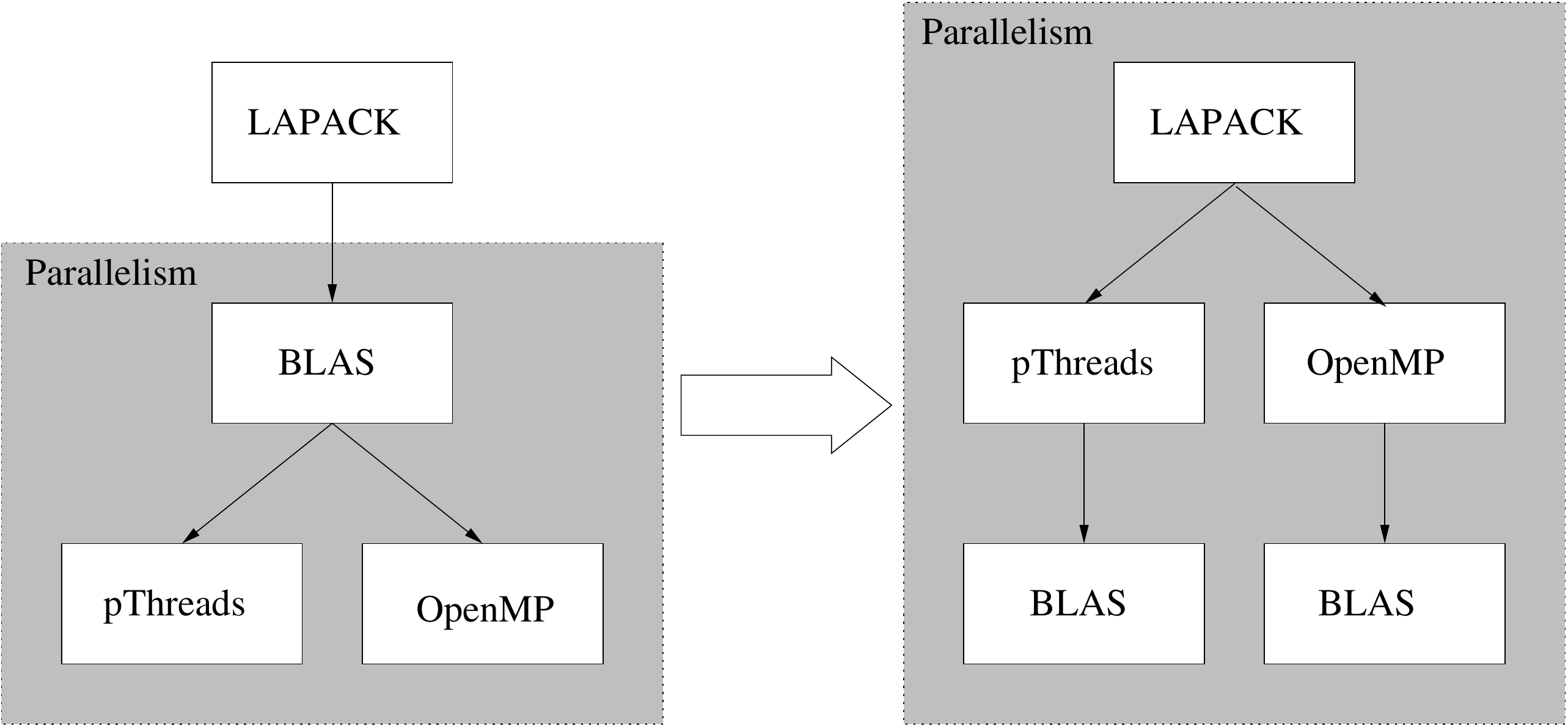}
  \caption{\label{fig:par_algos}Transition from sequential algorithms
    that rely on parallel BLAS to parallel algorithms.}
  \end{center}
\end{figure}

The LAPACK~\cite{lapack:99} and ScaLAPACK~\cite{scalapack:96} software
libraries represent a {\it de facto} standard for high performance
dense Linear Algebra computations and have been developed,
respectively, for shared-memory and distributed-memory
architectures. In both cases exploitation of parallelism comes from
the availability of parallel BLAS. In the LAPACK case, a number of
BLAS libraries can be used to take advantage of multiple processing
units on shared memory systems; for example, the freely distributed
ATLAS~\cite{ATLAS} and GotoBLAS~\cite{gotoblas} or other vendor BLAS
like Intel MKL~\cite{mkl} are popular choices. These parallel BLAS
libraries use common techniques for shared memory parallelization like
pThreads~\cite{mueller93pthreads} or OpenMP~\cite{10.1109/99.660313}.
This is represented in Figure~\ref{fig:par_algos} ({\it left}).

In the ScaLAPACK case, parallelism is exploited by
PBLAS~\cite{666023} which is a parallel BLAS implementation that uses
the Message Passing Interface~\cite{mpif1994} (MPI) for communications
on a distributed memory system.
Substantially, both LAPACK and ScaLAPACK implement sequential
algorithms that rely on parallel building blocks (i.e., the BLAS
operations).
As multicore systems require finer granularity and higher
asynchronicity, considerable advantages may be obtained by
reformulating old
algorithms or developing new algorithms in a way that their implementation can be
easily mapped on these new architectures. This transition is shown in
Figure~\ref{fig:par_algos}. An approach along these lines has already
been proposed in~\cite{Kurzak:2006:ILA,para06,1248397} where operations in
the standard LAPACK algorithms for some common factorizations were broken into sequences of smaller
tasks in order to achieve finer granularity and higher
flexibility in the scheduling of tasks to cores. The importance of
fine granularity algorithms is also shown in~\cite{cell_chol}.

The rest of this document shows how this can be achieved for the QR
factorization.
Section~\ref{sec:lapack} describes the algorithm for block QR factorization used
in the LAPACK
library; Section~\ref{sec:tiled} describes the tiled QR factorization
that provides both fine granularity and high level of asynchronicity;
performance results for this algorithm are shown in
Section~\ref{sec:perf}. Finally future working directions are
illustrated in Section~\ref{sec:future}.

\section{Block QR Factorization}
\label{sec:lapack}

\subsection{Description of the block QR Factorization}

The QR factorization is a transformation that factorizes an $m \times
n$ matrix $A$ into its factors $Q$ and $R$ where $Q$ is a unitary
matrix of size $n \times n$ and $R$ is a triangular matrix of size $m
\times m$. This factorization is operated by applying $min(m,n)$ Householder
reflections to the matrix $A$. Since Householder reflections are
orthogonal transformations, this factorization is stable as opposed to
the LU one; however, stability comes at the price of a higher flop
count: QR requires $2n^2(m-n/3)$ as opposed to the $n^2(m-n/3)$ needed for
LU. A detailed discussion of the QR factorization can be found
in~\cite{golubvanloan,trefethenbau,stew:98}. LAPACK uses a particular
version of this algorithm which achieves higher performance on
architectures with memory hierarchies thanks to blocking. This
algorithm is based on accumulating a number of Householder
transformations in what is called a {\it panel factorization} which are,
then, applied all at once by means of high performance Level 3 BLAS
operations. The technique used to accumulate Householder
transformation was introduced in~\cite{64889} and is known as the compact WY technique.

The LAPACK subroutine that performs the QR factorization is called
\texttt{DGEQRF} and is explained below.
Consider a matrix $A$ of size $m \times n$ that can be represented as
\begin{displaymath}
  A=\left(\begin{array}{cc}
    A_{11} & A_{12}  \\
    A_{21} & A_{22}  \\
  \end{array}\right)
\end{displaymath}
where $A_{11}$ is of size $b \times b$, $A_{12}$ of size $b \times
(n-b)$, $A_{21}$ of size $(m-b) \times b$ and $A_{22}$ of size $(m-b)
\times (n-b)$.

The LAPACK algorithm for QR factorization can be described as a
sequence of steps where, at each step, the transformation in
Equation~\eqref{eq:lap_qr0} is performed.

\begin{equation}
  \label{eq:lap_qr0}
  A=\left(\begin{array}{cc}
    A_{11} & A_{12}  \\
    A_{21} & A_{22}  \\
  \end{array}\right) \Longrightarrow
  \left(\begin{array}{c}
      V_{11} \\
      V_{21}
\end{array}\right), \left(\begin{array}{cc}

    R_{11}   & R_{12}        \\
      0     &\tilde{A}_{22} \\
  \end{array}\right)
\end{equation}

The transformation in Equation~\eqref{eq:lap_qr0} is obtained in two
steps:
\begin{enumerate}
\item {\bf Panel Factorization.} At this step a QR transformation of
  the panel $(A_{*1})$ is performed as in Equation~\eqref{eq:lap_qr1}.
  \begin{equation}
    \label{eq:lap_qr1}
    \left(\begin{array}{c}
        A_{11} \\
        A_{21}
      \end{array}\right) \Longrightarrow
    \left(\begin{array}{c}
        V_{11} \\
        V_{21}
      \end{array}\right), (T_{11}), (R_{11})
  \end{equation}
  This operation produces $b$ Householder reflectors $(V_{*,1})$ and
  an upper triangular matrix $R_{11}$ of size $b \times b$, which is a
  portion of the final $R$ factor, by means of the
  \texttt{DGEQR2} LAPACK subroutine; also, at this step, a triangular matrix $T_{11}$
  of size $b \times b$  by means of the \texttt{DLARFT} LAPACK
  subroutine\footnote{for the meaning of the matrix $T_{11}$ please
    refer to~\cite{64889}}. Please note that $V_{11}$ is a unit lower
  triangular matrix of size $b \times b$. The arrays $V_{*1}$ and
  $R_{11}$ do not need extra space to be stored since they overwrite
  $A_{*1}$. A temporary workspace is needed to store $T_{11}$.

\item {\bf Trailing submatrix update.} At this step, the
  transformation that was computed in the panel factorization is
  applied to the rest of the matrix, also called {\it trailing
    submatrix} as shown in Equation~\eqref{eq:lap_qr2}.
  \begin{equation}
    \label{eq:lap_qr2}
    \left(\begin{array}{c}
        R_{12}       \\
        \tilde{A}_{22}\\
      \end{array}\right) =
    \left(\begin{array}{c}
        \\
        \\
      \end{array}\right.  I - \left(\begin{array}{c}
        V_{11} \\
        V_{21}
      \end{array}\right) \begin{array}{c}
      \cdot (T_{11}) \cdot\\
      \\
    \end{array}   \begin{array}{cc}
      (V^T_{11}& V^T_{21}) \\
      &          \\
    \end{array}\left.\begin{array}{c}
        \\
        \\
      \end{array}\right)  \left(\begin{array}{c}
        A_{12}  \\
        A_{22}  \\
      \end{array}\right)
  \end{equation}
  This operation, performed by means of the \texttt{DLARFB} LAPACK
  subroutine, produces a portion $R_{12}$ of the final $R$ factor of
  size $b \times (n-b)$ and the matrix $\tilde{A}_{22}$.
\end{enumerate}
The QR factorization is continued by applying the
transformation~\eqref{eq:lap_qr0} to the submatrix $\tilde{A}_{22}$
and, then, iteratively, until the end of the matrix $A$ is reached. The
value of $b \ll m,n$ (the so called {\it block size}) is set by default to
32 in LAPACK-3.1.1 but different values may be more appropriate, and provide higher
performance, depending on the architecture characteristics.

\subsection{Scalability limits of the LAPACK implementation}
The LAPACK algorithm for QR factorization can use any flavor of parallel BLAS to exploit
parallelism on a multicore, shared-memory architecture.
This approach, however, has a number of limitations due to the nature
of the transformation in Equation~\eqref{eq:lap_qr1}, i.e., the panel
factorization. Both the \texttt{DGEQR2} and the \texttt{DLARFT} are
rich in Level 2 BLAS operations that cannot be efficiently parallelized on
currently available shared memory machines. To understand this, it is
important to note that Level 2 BLAS operations can be, generally speaking,
defined as all those operations where $O(n^2)$ floating-point operations
are performed on $O(n^2)$ floating-point data; thus, the speed of Level 2 BLAS
computations is limited by the speed at which the memory bus can feed
the cores. On current multicores architectures, there is a vast
disproportion between the bus bandwidth and the speed of the
cores. For example the Intel Clovertown processor is equipped with four cores each
capable of a double precision peak performance of 10.64 GFlop/s (that is to say a peak of
42.56 GFlop/s for four cores) while the bus bandwidth peak is 10.64 GB/s which provides 1.33
GWords/s (a word being a 64 bit double precision number).
As a result, since one core is largely enough to saturate
the bus, using two or more cores does not provide any significant
benefit. The LAPACK algorithm for QR factorization is, thus,
characterized by the presence of a sequential operation (i.e., the
panel factorization) which represents a small fraction
of the total number of FLOPS performed ($\mathcal{O}(n^2)$ FLOPS
for a total of $\mathcal{O}(n^3)$ FLOPS) but
limits the scalability of the block QR factorization on a multicore
system when parallelism is only exploited at the level of the BLAS
routines. This approach will be referred to as the {\it   fork-join}
approach since the execution flow of the QR factorization would show a
sequence of sequential operations (i.e. the panel factorizations)
interleaved to parallel ones (i.e., the trailing submatrix updates).

\begin{table}[!h]
  \centering
  \begin{tabular}[!h]{p{0.5in}p{0.5in}p{0.5in}p{0.3in}p{0.5in}p{0.5in}}
        \hline
        & \multicolumn{2}{c}{Intel MKL-9.1}  &  & \multicolumn{2}{c}{GotoBLAS-1.15}\\
        \hline
     \# cores   &  \texttt{DGEQR2} & \texttt{DGEQRF} &  & \texttt{DGEQR2} & \texttt{DGEQRF}\\
        &  Gflop/s         &   Gflop/s       &  & Gflop/s         & Gflop/s        \\
        \hline
        \hline
      1 &  0.4106          &  2.9            &  &  0.4549        &  3.31       \\
      2 &  0.4105          &  4.95           &  &  0.4558        &  5.51       \\
      4 &  0.4105          &  8.79           &  &  0.4557        &  9.69       \\
      8 &  0.4109          &  14.3           &  &  0.4549        &  10.58       \\
     16 &  0.4103          &  16.89          &  &  0.4558        &  13.01       \\
        \hline
  \end{tabular}
  \caption{\label{tab:blas2}Scalability of the fork-join
    parallelization on a 8-way dual Opteron system (sixteen cores
    total).}
\end{table}

Table~\ref{tab:blas2} shows the scalability limits of the panel
factorization and how this affects the scalability of the whole QR
factorization on an 8-way dual-core AMD Opteron system with MKL-9.1 and GotoBLAS-1.15
parallel BLAS libraries.

In~\cite{Kurzak:2006:ILA,para06}, a solution to this scalability
problem is presented. The approach described in~\cite{Kurzak:2006:ILA,para06} consists of
breaking the trailing submatrix update into smaller tasks that operate
on a block-column (i.e., a set of $b$ contiguous columns where $b$ is
the block size). The algorithm can then be represented as a Directed
Acyclic Graph (DAG) where nodes represent tasks, either panel factorization
or update of a block-column, and edges represent dependencies among
them. The execution of the algorithm is performed by asynchronously
scheduling the tasks in a way that dependencies are not violated. This
asynchronous scheduling results in an out-of-order execution where
slow, sequential tasks are hidden behind parallel ones. This approach
can be described as a dynamic lookahead technique.
Even if this approach provides significant speedup, as shown
in~\cite{Kurzak:2006:ILA}, it is exposed to scalability problems. Due
to the relatively high granularity of the tasks, the scheduling of
tasks may have a limited flexibility and the parallel execution of the
algorithm may be affected by an unbalanced load.

The following sections describe the application of this  idea of
dynamic scheduling and out of order execution to an algorithm for QR
factorization where finer granularity of the operations and higher
flexibility for the scheduling can be achieved.

\section{Tiled QR Factorization}
\label{sec:tiled}

The idea of dynamic scheduling and out of order execution can be
applied to a class of algorithms that allow the parallelization of
common Linear Algebra operations.  Previous work in this direction
includes SYRK, CHOL, block LU, and block QR
\cite{Kurzak:2006:ILA,1248397,para06}. For those four
factorizations, no algorithmic change is needed however, while CHOL
and SYRK can be naturally ``tiled'', the algorithms for block LU and
block QR factorizations involved tall and skinny panel factorization
that represents the bottlenecks of the computation (see pervious
Section). In order to have a finer granularity in LU and QR, we need
to ``tile'' the operations. To do so we will need a major
algorithmic change in LU and QR.

The algorithmic change we propose is actually well-known and takes
its roots in updating factorizations~\cite{golubvanloan,stew:98}.
Using updating techniques to tile the algorithms have
first\footnote{to our knowledge} been proposed by Yip
~\cite{yip_ooc} for LU to improve the efficiency of out-of-core
solvers, and were recently reintroduced in~\cite{vdgooclu,1055534}
for LU and QR, once more in the out-of-core context. A similar idea
has also been proposed in~\cite{210517} for Hessenberg reduction in
the parallel distributed context.

The originality of this paper is to study this techniques in the multicore
context, where they enable us to schedule operations of very fine
granularity.


\subsection{A Fine-Grain Algorithm for QR Factorization}

The QR factorization by blocks will be constructed based on the
following four elementary operations:
\begin{description}
\item[\texttt{DGEQT2}.] This subroutine was developed to perform the
  unblocked factorization of a diagonal block $A_{kk}$ of size $b \times
  b$. This operation produces an upper triangular matrix $R_{kk}$, a
  unit lower triangular matrix $V_{kk}$ that contains $b$ Householder
  reflectors and an upper triangular matrix $T_{kk}$ as defined by the
  WY technique for accumulating the transformations. Note that both
  $R_{kk}$ and $V_{kk}$ can be written on the memory area that was used
  for $A_{kk}$ and, thus, no extra storage is needed for them. A
  temporary work space is needed to store $T_{kk}$.

  Thus, \texttt{DGEQT2}($A_{kk}$, $T_{kk}$) performs
  \begin{displaymath}
    A_{kk} \longleftarrow V_{kk}, R_{kk} \qquad T_{kk} \longleftarrow T_{kk}
  \end{displaymath}

\item[\texttt{DLARFB}.] This LAPACK subroutine will be used to apply
  the transformation $(V_{kk}, T_{kk})$ computed by subroutine
  \texttt{DGEQT2} to a block $A_{kj}$.

  Thus, \texttt{DLARFB}($A_{kj}$, $V_{kk}$, $T_{kk}$) performs
  \begin{displaymath}
    A_{kj} \longleftarrow (I-V_{kk}T_{kk}V^T_{kk})A_{kj}
  \end{displaymath}

\item[\texttt{DTSQT2}.] This subroutine was developed to perform the
  unblocked QR factorization of a matrix that is formed by coupling an
  upper triangular block $R_{kk}$ with a square block $A_{ik}$. This
  subroutine will return an upper triangular matrix $\tilde{R}_{kk}$
  which will overwrite $R_{kk}$ and $b$ Householder reflectors where $b$
  is the block size. Note
  that, since $R_{kk}$ is upper triangular, the resulting Householder
  reflectors can be represented as an identity block $I$ on top of a
  square block $V_{ik}$. For this reason no extra storage is needed for
  the Householder vectors  since the identity block need not be stored
  and $V_{ik}$ can overwrite $A_{ik}$. Also a matrix $T_{ik}$ is
  produced for which storage space has to be allocated.

  Then, \texttt{DTSQT2}($R_{kk}$, $A_{ik}$, $T_{ik}$) performs
  \begin{displaymath}
    \left(
      \begin{array}[!h]{c}
        R_{kk}\\
        A_{ik}
      \end{array}\right) \longleftarrow
    \left(\begin{array}[!h]{c}
        I\\
        V_{ik}
      \end{array}\right), \tilde{R}_{kk} \qquad T_{ik} \longleftarrow T_{ik}
  \end{displaymath}

\item[\texttt{DSSRFB}.] This subroutine was developed to apply the
  transformation computed by \texttt{DTSQT2} to a matrix formed
  coupling two square blocks $A_{kj}$ and $A_{ij}$.

  Thus, \texttt{DSSRF}($A_{kj}$, $A_{ij}$, $V_{ik}$, $T_{ik}$)
  performs
  \begin{displaymath}
    \left(\begin{array}{c}
      A_{kj}\\
      A_{ij}
    \end{array}\right) \longleftarrow
  \left(
    \begin{array}{c}
      \\
      \\
    \end{array}\right. I -
    \left(
      \begin{array}{c}
        I \\
        V_{ik}
      \end{array}\right)
    \begin{array}{c}
     \cdot (T_{ik}) \cdot\\
      \\
    \end{array}
    \begin{array}{cc}
      (I & V^T_{ik})\\
      \\
    \end{array}  \left.
    \begin{array}{c}
      \\
      \\
    \end{array}\right)\left(\begin{array}{c}
      A_{kj}\\
      A_{ij}
    \end{array}\right)
  \end{displaymath}
\end{description}

All of this elementary operations rely on BLAS subroutines to perform
internal computations.

Assuming a matrix $A$ of size $pb \times qb$
\begin{displaymath}
  \left(\begin{array}[!h]{cccc}
    A_{11}  & A_{12} & \dots  & A_{1q}  \\
    A_{21}  & A_{22} & \dots  & A_{2q}  \\
    \vdots &       & \ddots & \vdots \\
    A_{p1}  & A_{p2} & \dots  & A_{pq}
  \end{array}\right)
\end{displaymath}
where $b$ is the block size and each $A_{ij}$ is of size $b \times b$,
the QR factorization can be performed as in Algorithm~\ref{alg:blkqr}.

\begin{algorithm}[!h]
\caption{\label{alg:blkqr}The block algorithm for QR factorization.}
  \begin{algorithmic}[1]
    \FOR{$k=1, 2..., min(p,q)$}
    \STATE \texttt{DGEQT2}($A_{kk}$, $T_{kk}$);
    \FOR{$j=k+1, k+2, ..., q$}
    \STATE \texttt{DLARFB}($A_{kj}$, $V_{kk}$, $T_{kk}$);
    \ENDFOR
    \FOR{$i=k+1, k+1, ..., p$}
    \STATE \texttt{DTSQT2}($R_{kk}$, $A_{ik}$, $T_{ik}$);
    \FOR{$j=k+1, k+2, ..., q$}
    \STATE \texttt{DSSRFB}($A_{kj}$, $A_{ij}$, $V_{ik}$, $T_{ik}$);
    \ENDFOR
    \ENDFOR
    \ENDFOR
  \end{algorithmic}
\end{algorithm}

\begin{figure}[!h]
  \begin{center}
    \includegraphics[width=\textwidth]{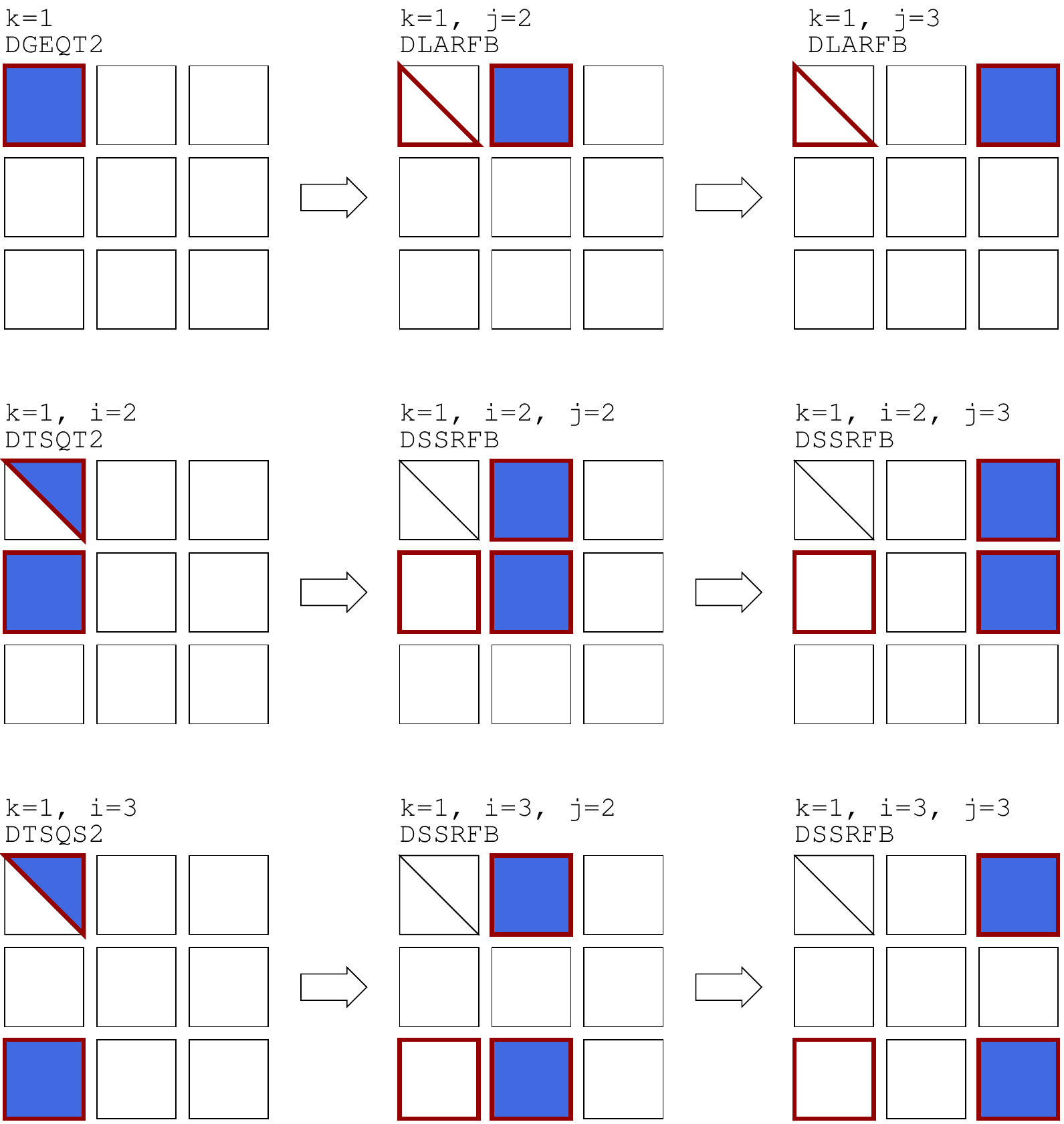}
  \caption{\label{fig:blk_alg}Graphical representation of one
    repetition of the outer loop in Algorithm~\ref{alg:blkqr} on a
    matrix with $p=q=3$. As expected the picture is very similar
    to the out-of-core algorithm presented in~\cite{1055534}.
}
  \end{center}
\end{figure}

Figure~\ref{fig:blk_alg} gives a graphical representation of one
repetition (with $k=1$) of the outer loop in Algorithm~\ref{alg:blkqr} with
$p=q=3$. The red, thick borders show what blocks in the matrix are being read
and the light blue fill shows what blocks are being written in a
step. The $T_{kk}$ matrices are not shown in this figure for clarity purposes.

\subsection{Operation count}

This section shows that Algorithm~\ref{alg:blkqr} has a
higher operation count than the LAPACK algorithm discussed in
Section~\ref{sec:lapack}. Performance results in
Section~\ref{sec:perf} will demonstrate that it is worth paying this
cost for the sake of scaling.
The operation count of the block algorithm for QR factorization can be derived starting
from the operation count of each elementary operation; assuming that
$b$ is the block size:
\begin{description}
\item[\texttt{DGEQT2}:] this operation is a standard non blocked QR
  factorization of a $b \times b$ matrix where, in addition, the
  $T_{kk}$ triangular matrix is computed. Thus, $4/3 b^3$ floating
  point operations are performed for the factorization plus $2/3b^3$
  for computing $T_{kk}$. This subroutine accounts for $2b^3$ floating
  point operations total.

\item[\texttt{DLARFB}:] since both $V_{kk}$ and $T_{kk}$ are
  triangular matrices, $3b^3$ floating point operations are done in
  this subroutine.
\item[\texttt{DTSQT2}:] taking advantage of the triangular structure of
  $R_{kk}$, the factorization can be computed in $2b^3$ floating point
  operations. The computation of the triangular $T_{ik}$ matrix can
  also be performed exploiting the structure of the Householder
  vectors built during the factorization (remember that the $b$
  reflectors can be represented as an identity block on top of a
  square full block). Since $4/3b^3$ are needed to compute $T_{ik}$,
  the \texttt{DTSQT2} accounts for $10/3b^3$ floating point operations.
\item[\texttt{DSSRFB}:] exploiting the structure of the Householder
  reflectors and of the $T_{ik}$ matrix computed in \texttt{DTSQT2},
  this subroutine needs $5b^3$ floating point operations.
\end{description}

For each repetition of the outer loop in Algorithm~\ref{alg:blkqr},
one \texttt{DGEQT2}, $q-k$ \texttt{DLARFB}, $p-k$ \texttt{DTSQT2} and
$(p-k)(q-k)$ \texttt{DSSRFB} are performed for a total of
$2b^3+3(q-k)b^3+10/3(p-k)b^3+5(p-k)(q-k)b^3$. Assuming, without loss
of generality, that $q<p$ and integrating this quantity over all the
$q$ repetitions of the outer loop in Algorithm~\ref{alg:blkqr}, the
total operation count for the QR factorization is
\begin{equation}
  \label{eq:opcount}
  \begin{array}{l}
  \sum_{k=1}^{q}(2b^3+3(q-k)b^3+\frac{10}{3}(p-k)b^3+5(p-k)(q-k)b^3)\\
  \\
  \simeq \frac{5}{2}q^2(p-\frac{q}{3})b^3 \\
  \\
  = \frac{5}{2}n^2(m-\frac{n}{3}).
  \end{array}
\end{equation}

Equation~\eqref{eq:opcount} shows that the block algorithm for QR
factorization needs 25\% more floating point operations than the
standard LAPACK algorithm.

Note that if we have used the nonblocked version of the Householder application
(DLARF instead of DLARFB) then the number of FLOPS for the tiled algorithm
and the block algorithm would have been exactly the same.

An easy twist to reduce the 25\% overhead of the block version is to use nonsquare
blocks. For example, using $2b\times b$ blocks would reduce the overhead at 12.5\%.

\subsection{Graph driven asynchronous execution}

\begin{figure}[!h]
  \begin{center}
    \includegraphics[width=\textwidth]{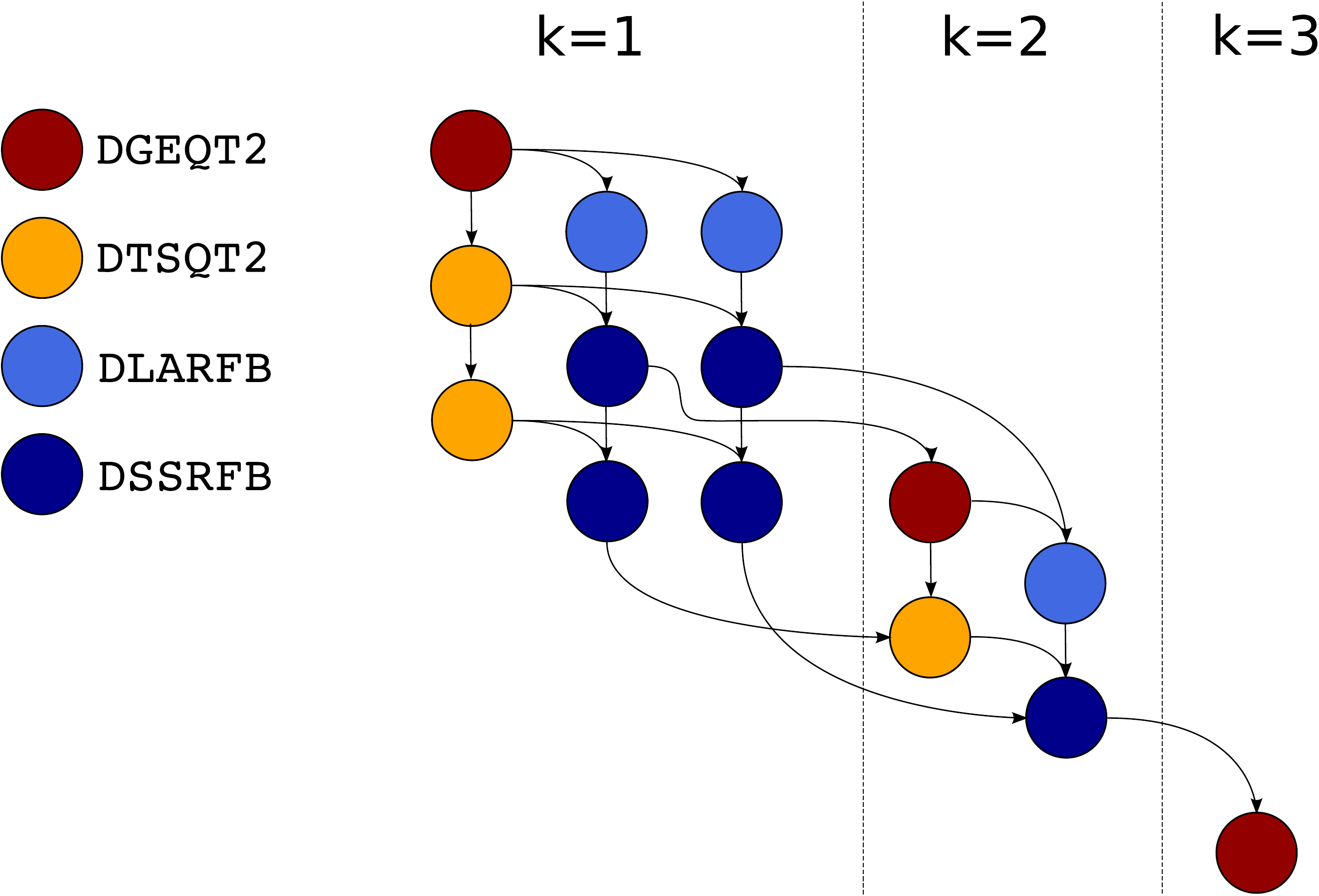}
  \caption{\label{fig:qr_dag}The dependency graph of
    Algorithm~\ref{alg:blkqr} on a matrix with $p=q=3$.}
  \end{center}
\end{figure}

Following the approach presented in~\cite{Kurzak:2006:ILA,para06},
Algorithm~\ref{alg:blkqr} can be represented as a Directed Acyclic
Graph (DAG) where nodes are elementary tasks that operate on
$b \times b$ blocks and where edges represent the dependencies among
them. Figure~\ref{fig:qr_dag} show the DAG when
Algorithm~\ref{alg:blkqr} is executed on a matrix with $p=q=3$. Note
that the DAG has a recursive structure and, thus, if $p_1 \ge p_2$ and
$q_1 \ge q_2$ then the DAG for a matrix of size $p_2 \times q_2$ is
a subgraph of the DAG for a matrix of size $p_1 \times q_1$. This
property also holds for most of the algorithms in LAPACK.

Once the DAG is known, the tasks can be scheduled asynchronously and
independently as long as the dependencies are not violated. A critical
path can be identified in the DAG as the path that connects all the
nodes that have the higher number of outgoing edges. Based on this
observation, a scheduling policy can be used, where higher priority is
assigned to those nodes that lie on the critical path. Clearly, in the
case of our block algorithm for QR factorization, the nodes associated to the
\texttt{DGEQT2} subroutine have the highest priority and then three other
priority levels can be defined for \texttt{DTSQT2}, \texttt{DLARFB} and
\texttt{DSSRFB} in descending order.

\begin{figure}[!h]
  \begin{center}
    \includegraphics[width=\textwidth]{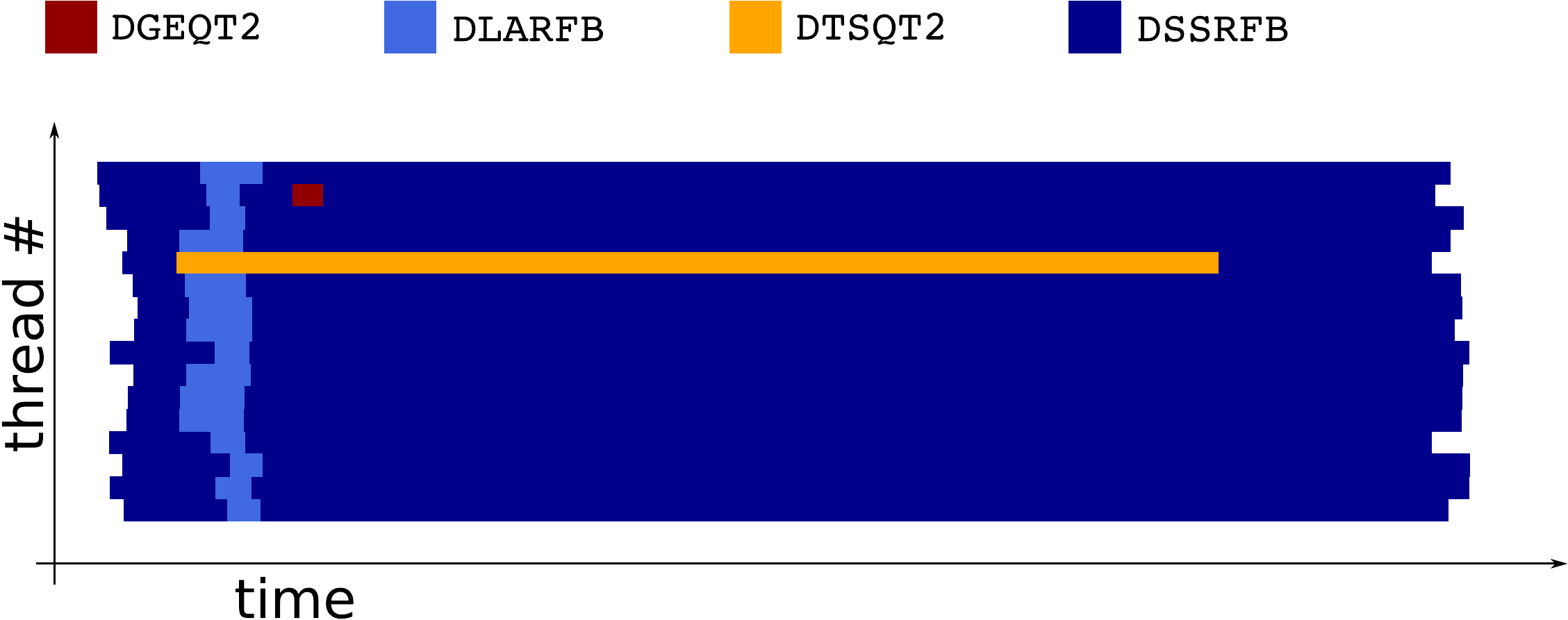}
  \caption{\label{fig:flow}The execution flow for dynamic scheduling,
    out of order execution of Algorithm~\ref{alg:blkqr}.}
  \end{center}
\end{figure}

This
dynamic scheduling results in an out of order execution where idle
time is almost completely eliminated since only very loose
synchronization is required between the threads. Figure~\ref{fig:flow}
shows part of the execution flow of Algorithm~\ref{alg:blkqr} on a
16-cores machine (8-way Dual Opteron) when tasks are dynamically
scheduled based on dependencies in the DAG. Each line in the execution
flow shows which tasks are performed by one of the threads involved in the
factorization.

Figure~\ref{fig:flow} shows that all the idle times, which represent
the major scalability limit of the fork-join approach, can be
 removed thanks to the very low synchronization requirements
of the graph driven execution. The graph driven execution also
provides some degree of adaptivity since tasks are scheduled to
threads depending on the availability of execution units.

\subsection{Block Data Layout}

The major limitation of performing very fine grain computations, is
that the BLAS library generally have very poor performance on small
blocks. This situation can be considerably improved by storing
matrices in Block Data Layout (BDL) instead of the Column Major Format
that is the standard storage format for FORTRAN arrays.

\begin{figure}[!h]
  \begin{center}
    \includegraphics[width=0.7\textwidth]{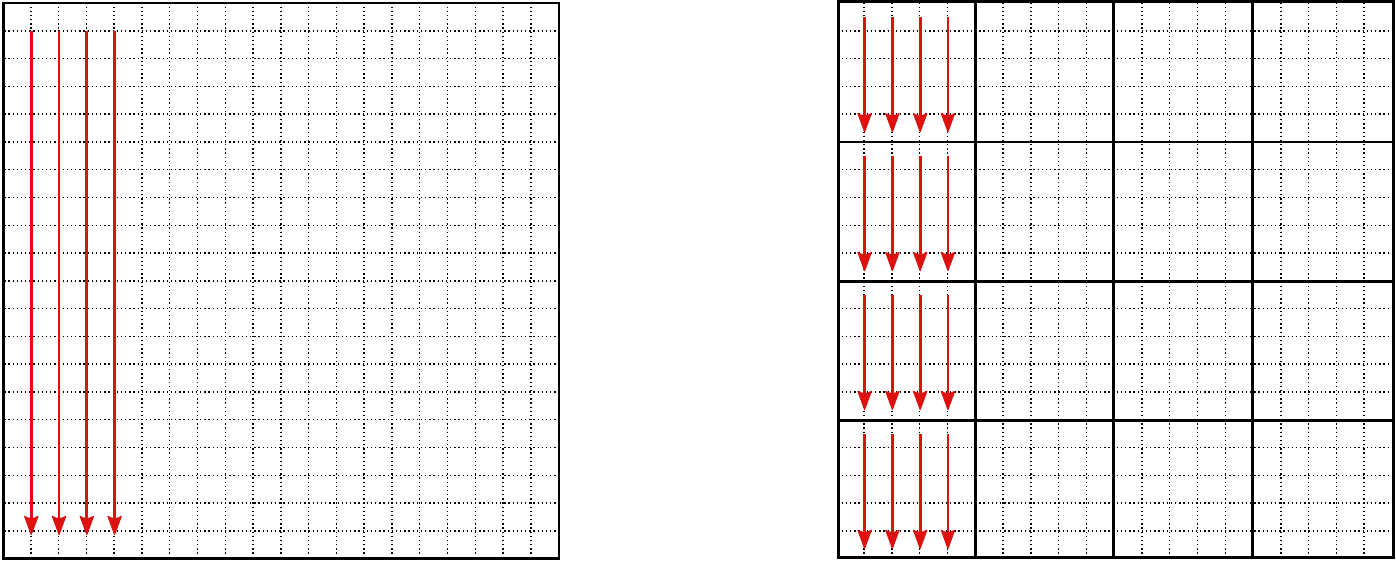}
  \caption{\label{fig:bdl}A comparison of Column Major storage format
    ({\it left}) and Block Data Layout ({\it right}).}
  \end{center}
\end{figure}

Figure~\ref{fig:bdl} compares Column Major Format ({\it left}) and
Block Data Layout ({\it right}). In BDL a matrix is split into blocks
and each block is stored into contiguous memory locations. Each block
is stored in Column Major Format and blocks are stored in Column Major
Format with respect to each other. As a result the access pattern to
memory is more regular and BLAS performance is considerably improved.
The benefits of BDL have been extensively studied in the past, for
example in~\cite{670985}, and recent studies like~\cite{1248397} demonstrate how
fine-granularity parallel algorithms can benefit from BDL.

\section{Performance Results}
\label{sec:perf}

The performance of the tiled QR factorization with dynamic scheduling of
tasks has been measured on the systems listed in
Table~\ref{tab:archs} and compared to the performance of the fork-join
approach, i.e., the standard algorithm for block QR factorization of LAPACK associated with
multithreaded BLAS.

\begin{table}[!h]
  \centering
  \begin{tabular}{p{1in}p{1.8in}p{1.8in}}
\hline
                 &  8-way dual Opteron          & 2-way quad Clovertown  \\
\hline
Architecture     &  Dual-Core AMD               & Intel\textregistered  Xeon\textregistered  CPU    \\
                 &  Opteron\texttrademark  8214            &           X5355         \\
Clock speed      &  2.2 GHz                     &   2.66 GHz              \\
\# cores         &  $8\times2=16$               &   $2\times4=8$          \\
Peak performance &  70.4  Gflop/s               &   85.12 Gflop/s         \\
Memory           &  62 GB                       &   16 GB                 \\
Compiler suite   &  Intel 9.1                   &   Intel 9.1             \\
BLAS libraries   &  GotoBLAS-1.15               &   GotoBLAS-1.15         \\
                 &  MKL-9.1                     &   MKL-9.1               \\
\hline
  \end{tabular}
  \caption{\label{tab:archs}Details of the systems used for the following performance results.}
\end{table}

\begin{figure*}
\begin{minipage}[tl]{0.5\textwidth}
\begin{center}
{\includegraphics[width=1\textwidth]{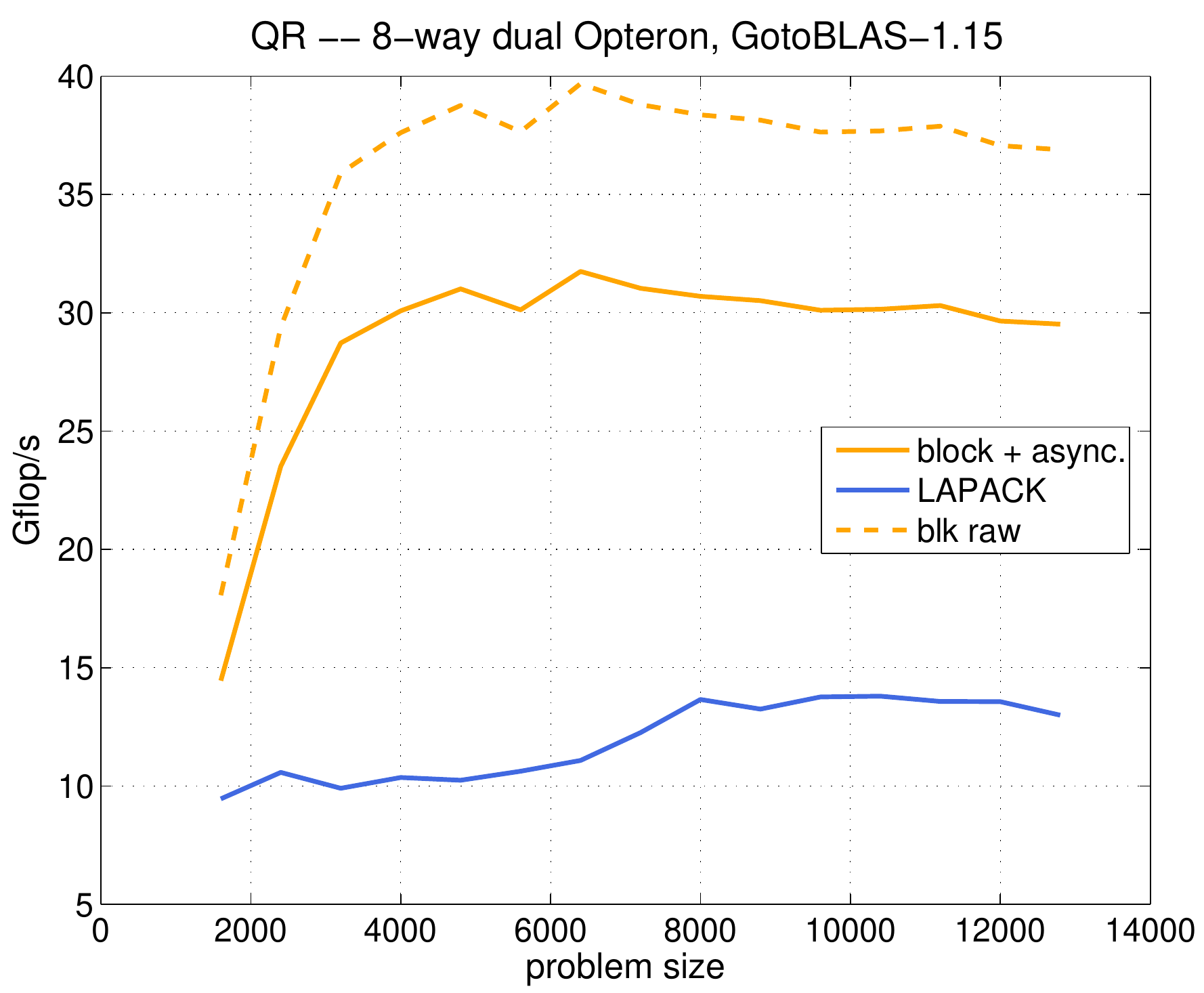}}
\end{center}
\end{minipage}
\hspace{0.25cm}
\begin{minipage}[tr]{0.5\textwidth}
\begin{center}
{\includegraphics[width=1\textwidth]{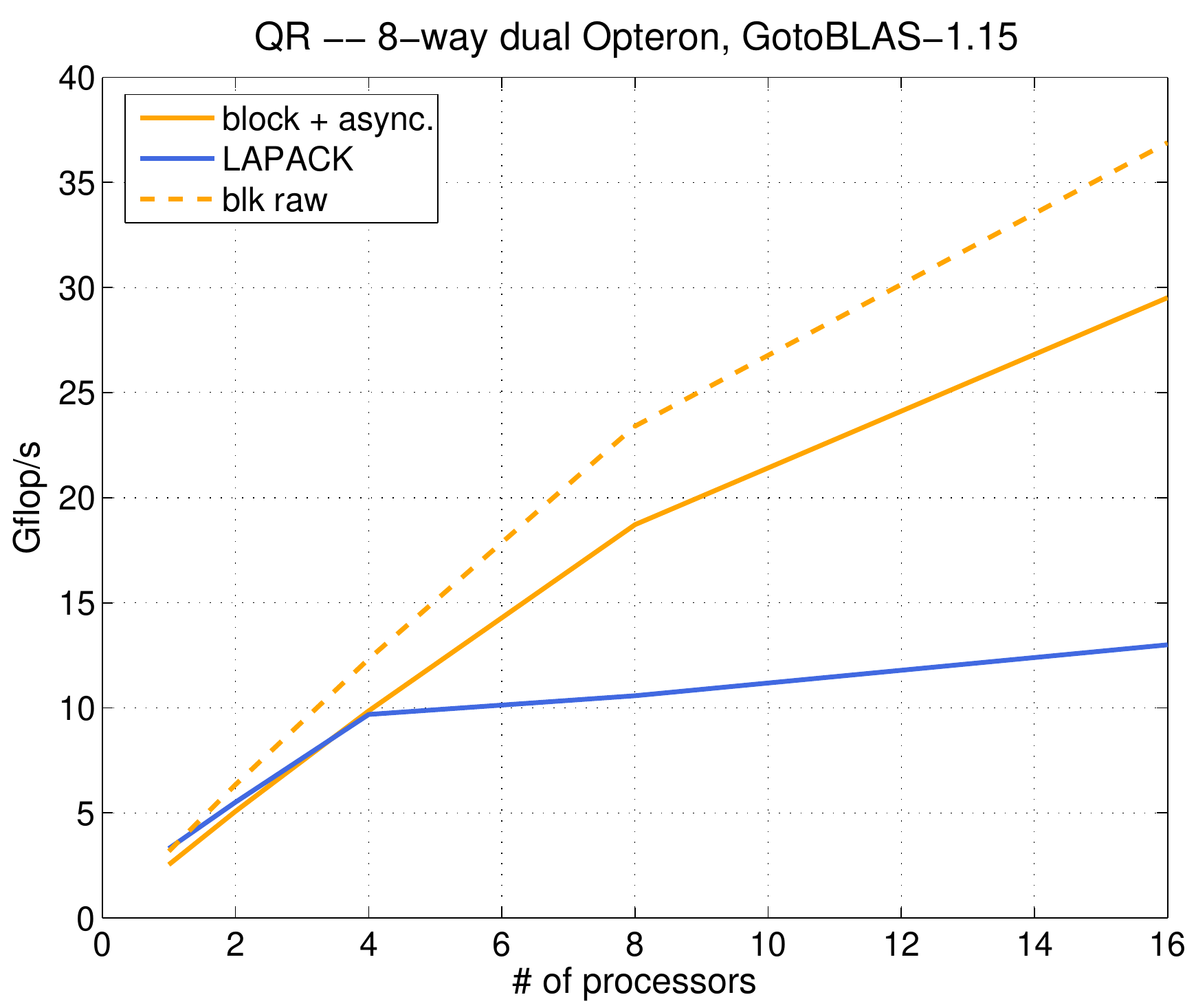}}
\end{center}
\end{minipage}
\caption{\label{fig:goto_opt}Comparison between the performance of the
block algorithm with dynamic scheduling using GotoBLAS-1.15 on an
8-way dual Opteron system. The dashed curve reports the
raw performance of the block algorithm with dynamic scheduling, i.e.,
the performance as computed with the true operation count in Equation~\eqref{eq:opcount}.}
\end{figure*}

\begin{figure*}
\begin{minipage}[tl]{0.5\textwidth}
\begin{center}
{\includegraphics[width=1\textwidth]{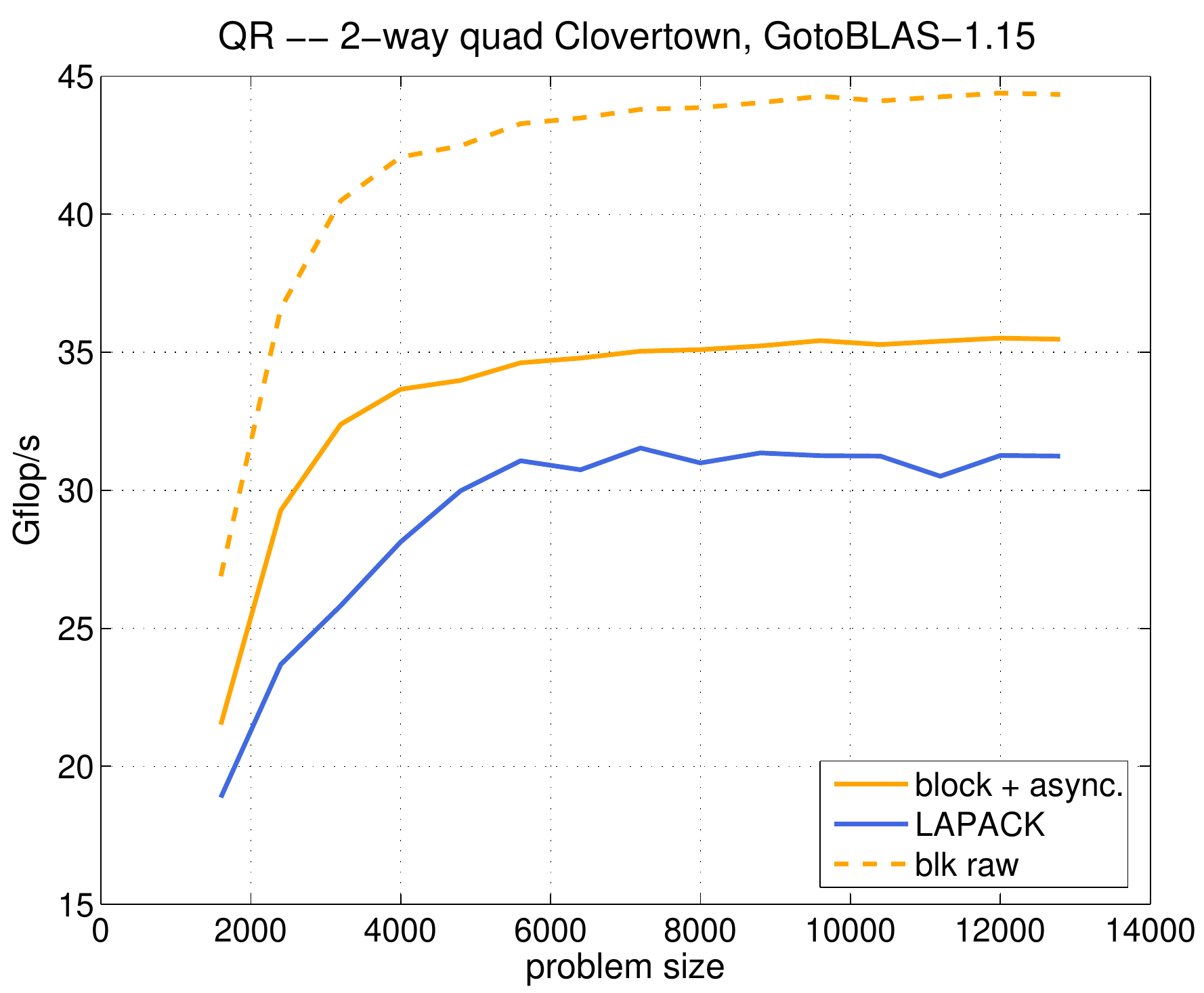}}
\end{center}
\end{minipage}
\hspace{0.25cm}
\begin{minipage}[tr]{0.5\textwidth}
\begin{center}
{\includegraphics[width=1\textwidth]{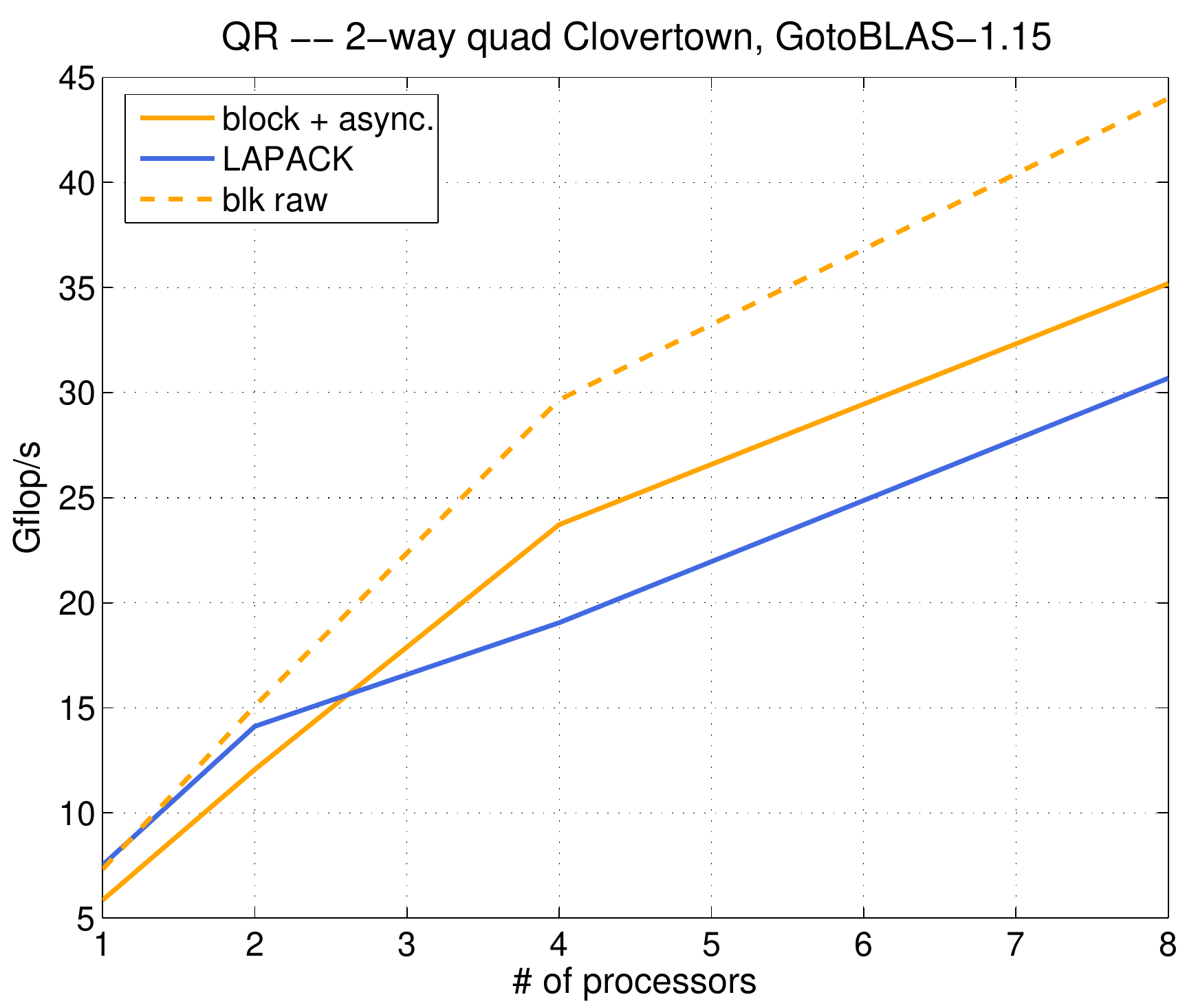}}
\end{center}
\end{minipage}
\caption{\label{fig:goto_clov}Comparison between the performance of the
block algorithm with dynamic scheduling using GotoBLAS-1.15 on an
2-way quad Clovertown system. The dashed curve reports the
raw performance of the block algorithm with dynamic scheduling, i.e.,
the performance as computed with the true operation count in Equation~\eqref{eq:opcount}.}
\end{figure*}

\begin{figure*}
\begin{minipage}[tl]{0.5\textwidth}
\begin{center}
{\includegraphics[width=1\textwidth]{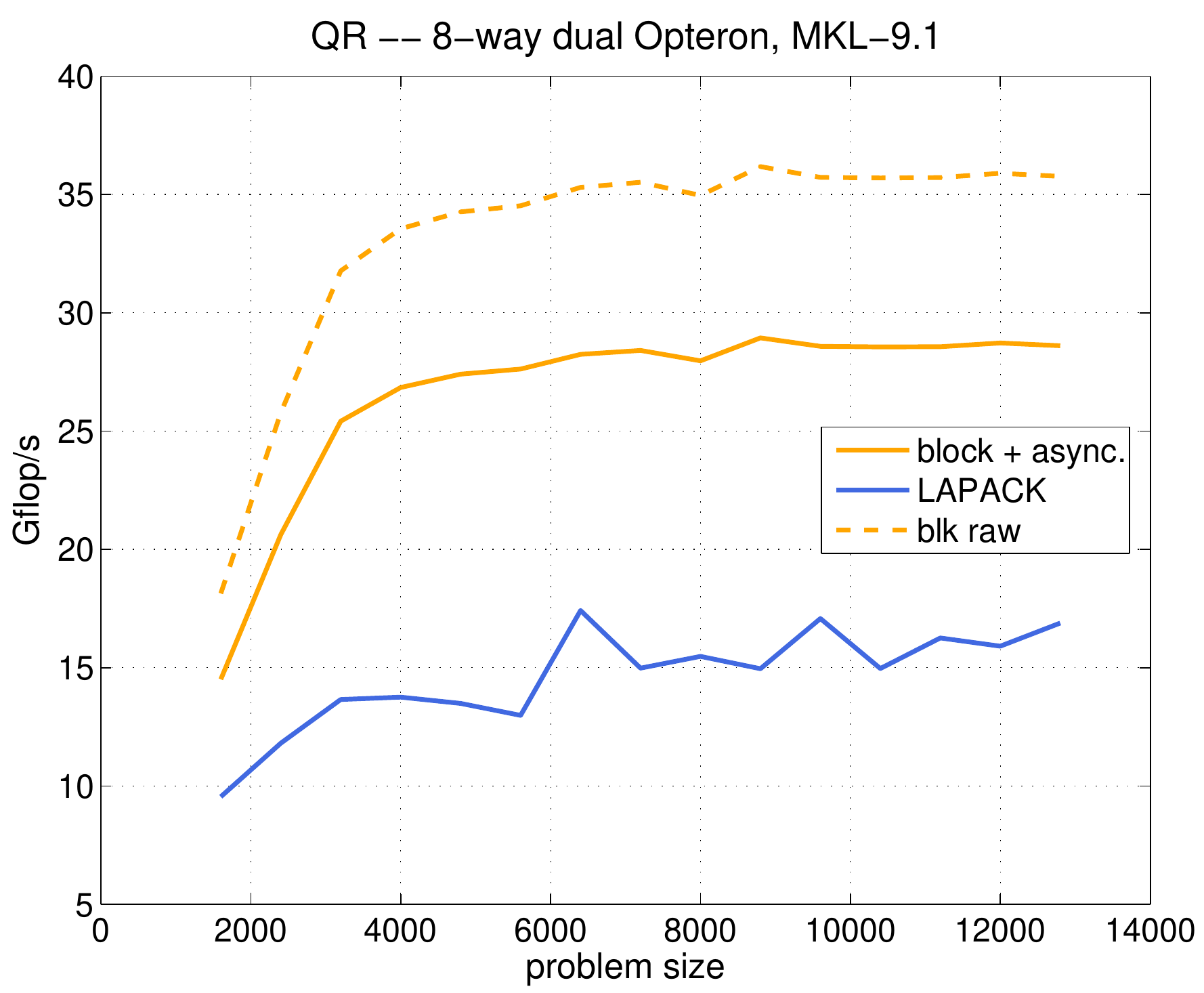}}
\end{center}
\end{minipage}
\hspace{0.25cm}
\begin{minipage}[tr]{0.5\textwidth}
\begin{center}
{\includegraphics[width=1\textwidth]{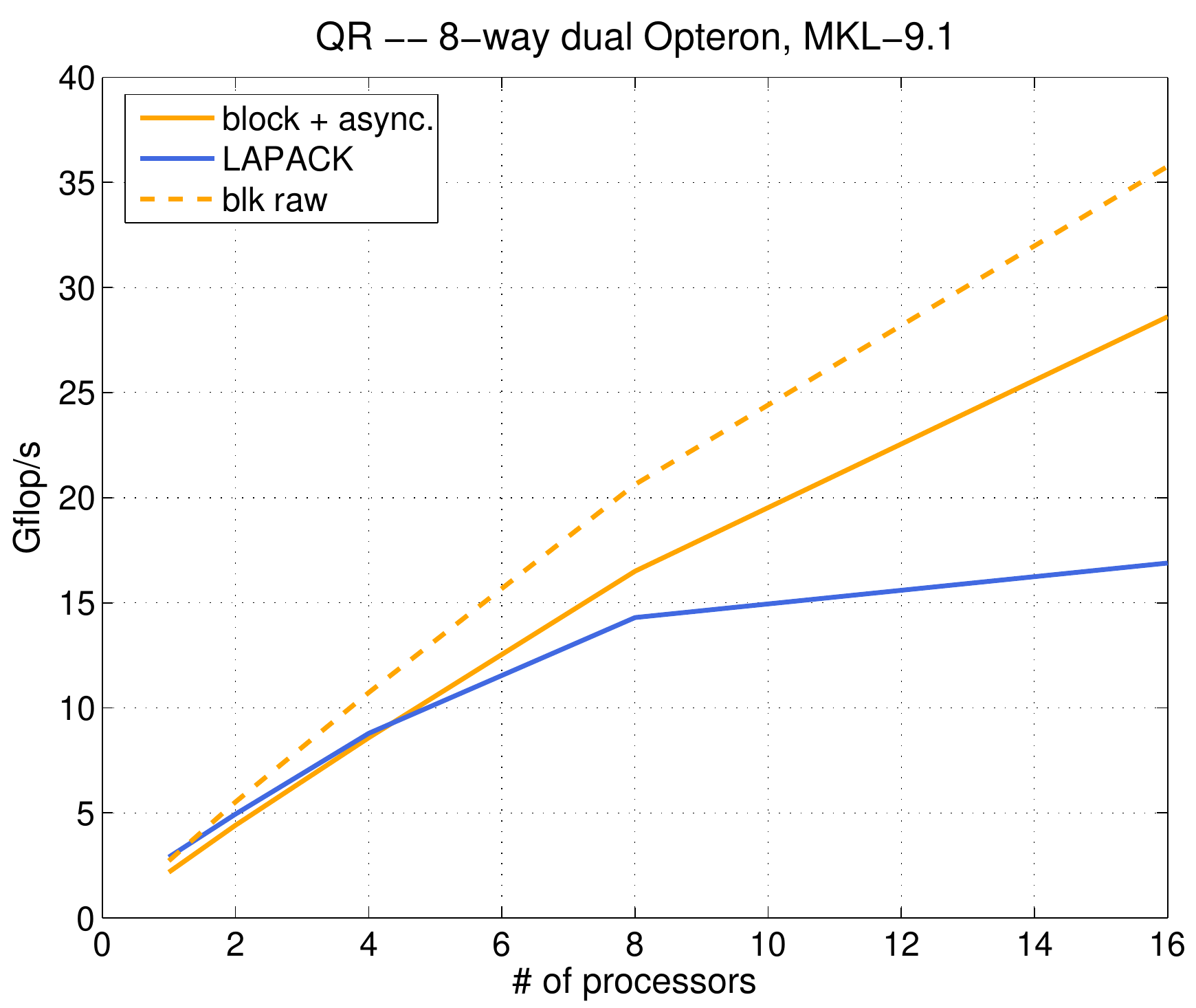}}
\end{center}
\end{minipage}
\caption{\label{fig:mkl_opt}Comparison between the performance of the
block algorithm with dynamic scheduling using MKL-9.1 on an
8-way dual Opteron system. The dashed curve reports the
raw performance of the block algorithm with dynamic scheduling, i.e.,
the performance as computed with the true operation count in Equation~\eqref{eq:opcount}.}
\end{figure*}

\begin{figure*}
\begin{minipage}[tl]{0.5\textwidth}
\begin{center}
{\includegraphics[width=1\textwidth]{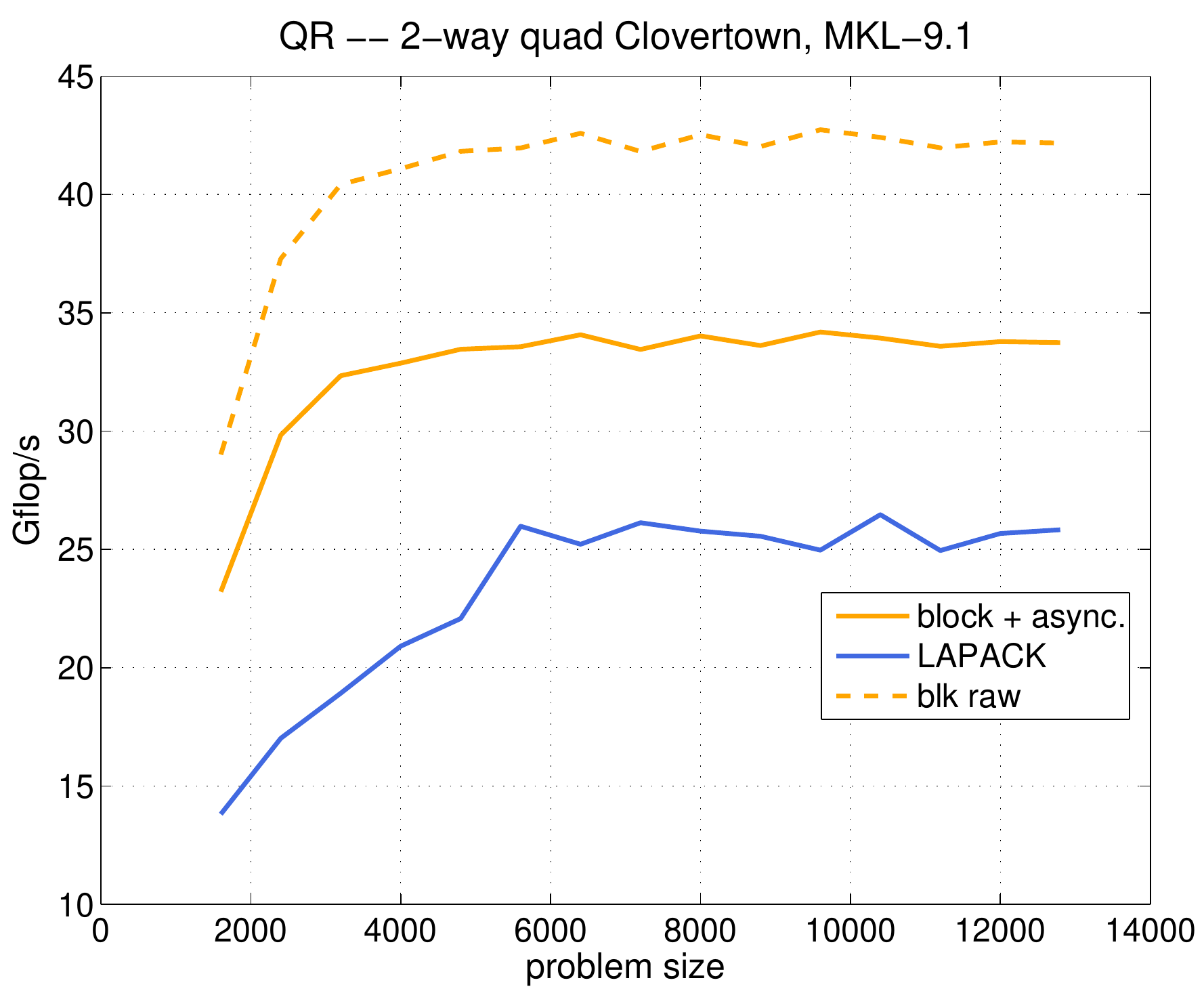}}
\end{center}
\end{minipage}
\hspace{0.25cm}
\begin{minipage}[tr]{0.5\textwidth}
\begin{center}
{\includegraphics[width=1\textwidth]{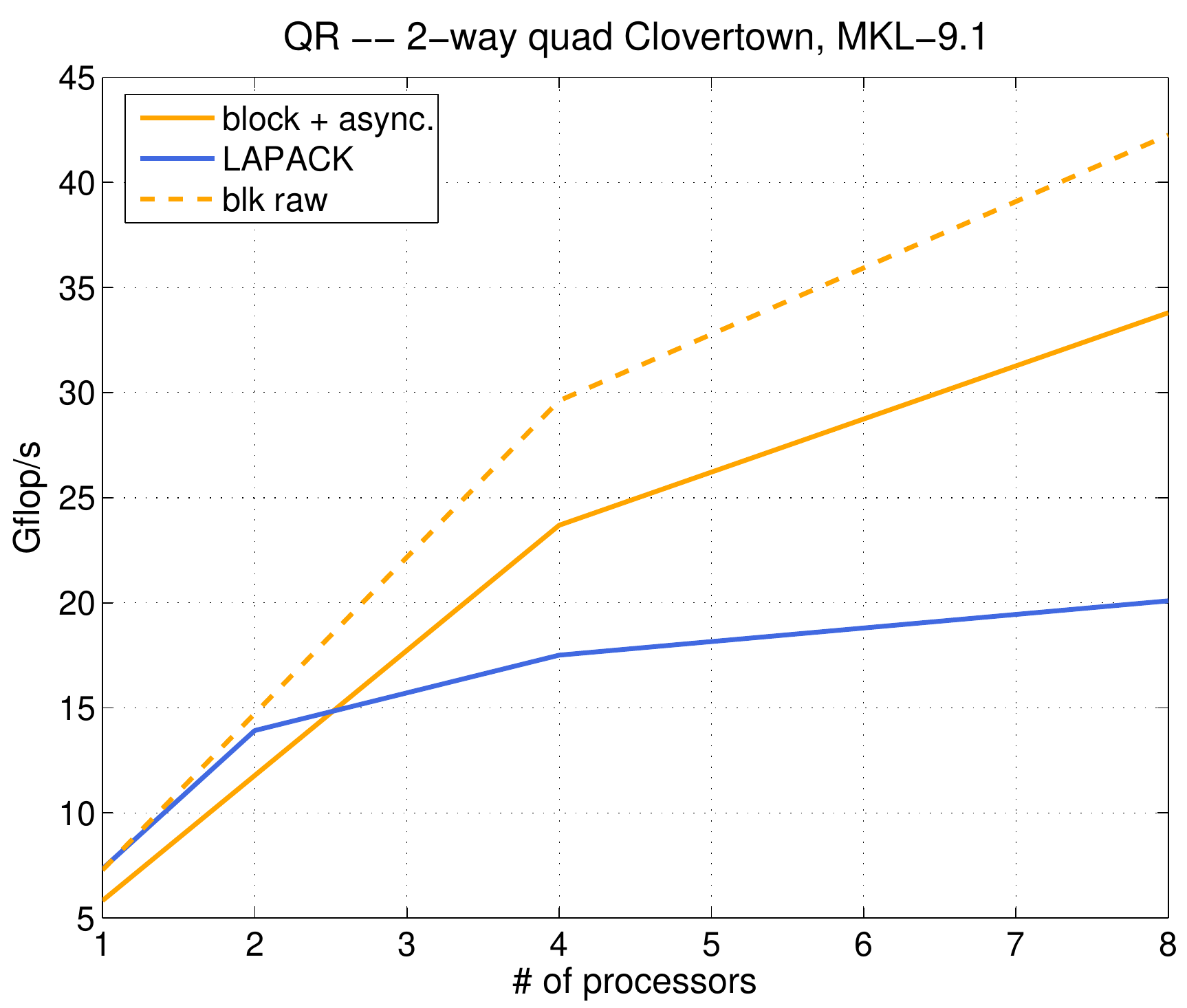}}
\end{center}
\end{minipage}
\caption{\label{fig:mkl_clov}Comparison between the performance of the
block algorithm with dynamic scheduling using MKL-9.1 on an
2-way quad Clovertown system. The dashed curve reports the
raw performance of the block algorithm with dynamic scheduling, i.e.,
the performance as computed with the true operation count in Equation~\eqref{eq:opcount}.}
\end{figure*}

Figures~\ref{fig:goto_opt},~\ref{fig:goto_clov},~\ref{fig:mkl_opt},~\ref{fig:mkl_clov}
report the performance of the QR factorization for both the block algorithm with dynamic
scheduling and the LAPACK algorithm with multithreaded BLAS. A block
size of 200 has been used for the block algorithm while the block size
for the LAPACK algorithm\footnote{the block size in the LAPACK
  algorithm sets the width of the panel.} has been tuned in order to
achieve the best performance for all the combinations of architecture
and BLAS library.

In each graph, two curves are reported for the
block algorithm with dynamic scheduling; the solid curve shows its
relative performance when the operation count is assumed equal to the one of the
LAPACK algorithm reported in Section~\ref{sec:lapack} while the dashed
curve shows its ``raw'' performance, i.e. the actual flop rate computed with
the exact operation count for this algorithm (given in Equation~\eqref{eq:opcount}).
As already mentioned,
the ``raw performance'' (dashed curve) is 25\% higher than the
relative performance (solid curve).

The graphs on the left part of each figure show the performance
measured using the maximum number of cores available on each system
with respect to the problem size. The graphs on the right part of each
figure show the weak scalability, i.e. the flop rates versus the
number of cores when the local problem size is kept constant
(nloc=5,000) as the number of cores increases.

Figures~\ref{fig:goto_opt},~\ref{fig:goto_clov},~\ref{fig:mkl_opt},~\ref{fig:mkl_clov}
show that, despite the higher operation count, the block algorithm
with dynamic scheduling is capable of completing the QR
factorization in less time than the LAPACK algorithm when the
parallelism degree is high enough that the benefits of the
asynchronous execution overcome the penalty of the extra flops. For
lower numbers of cores, in fact, the fork-join approach has a good
scalability and completes the QR factorization in less time than the
block algorithm because of the lower flop count. Note that the
actual execution rate of the block algorithm for QR factorization with dynamic
scheduling (i.e., the dashed curves) is always higher than that of the
LAPACK algorithm with multithreaded BLAS even for low numbers of cores.
The actual performance of the block algorithm, even if considerably
higher than that of the fork-join one, is still far from the peak
performance of the systems used for the measures. This is mostly due
to two factors. First the nature of the BLAS operations involved;
the \texttt{DGEQR2} and the \texttt{DLARFT} in the LAPACK algorithm
and the \texttt{DGEQT2} and \texttt{DTSQT2} in the block algorithm are
based on Level 2 BLAS operations that, being memory bound, represent a limit
for performance. Second, the performance of BLAS routines on small
size blocks. The block size used in the experiments reported above is
200; this block size represents a good compromise between flexibility
of the scheduler and performance of the BLAS operations but it is far
from being ideal. Such a block size, in fact, does not allow a good
task scheduling for smaller size problems and still the performance of
BLAS operations is far from what can be achieved for bigger size blocks.

\section{Conclusion}
\label{sec:future}

By adapting known algorithms for updating the QR factorization of a
matrix, we have derived an implementation scheme of the QR
factorization for multicore architectures based on dynamic
scheduling and block data layout. Although the proposed algorithm is
performing 25\% more FLOPS than the regular algorithm, the gain in
flexibility allows an efficient dynamic scheduling which enables the
algorithm to scale almost perfectly when the number of cores
increases.

We note that the 25\% overhead can be easily reduced by using nonsquare blocks.
For example, using $2b\times b$ blocks, the overhead reduces to 12.5\%.

While this paper only addresses the QR factorization, it is straightforward to
derive with the same ideas the two important computational routines that consists
in applying the Q factor to a set of vectors (see DORMQR in LAPACK) and constructing
the Q-factor (see DORGQR in LAPACK).

\noindent
The ideas behind this work can be extended in many directions:
\begin{description}
\item[Implement other linear algebra operations.] The LU factorization
  can be performed with an algorithm that is analogous to the QR one
  described in Section~\ref{sec:tiled}. This algorithm has been
  discussed in~\cite{yip_ooc,vdgooclu} as a way of improving the out-of-core LU factorization.
  Even though the only difference between the
  block algorithms for the LU and QR factorizations is in the elementary
  operations, in the LU case the cost of block algorithm is 50\%
  higher than the LAPACK algorithm. For this reason, the benefits of
  the improved scalability may be visible only at very high processor
  counts or may not be visible at all. Techniques must be investigated
  to eliminate or reduce the extra cost.

  The same to blocking may also be applied to other two sided
  transformations like Hessenberg reduction, Tridiagonalization and
  Bidiagonalization. In these transformations, Level 2 BLAS operations are
  predominant and panel reductions account for almost 50\% of the time
  of a sequential execution. Breaking the panel into smaller tasks
  that can be executed in parallel with other tasks may yield
  considerable performance improvements.

\item[Enforcing data locality.] The results proposed in~\cite{1248397}
  show that enforcing data locality and CPU affinity may provide
  considerable benefits. It must be noted that the improvements that
  can be expected on non-multicore SMPs are higher than on currently
  available multicore systems and this is due to the fact that
  on multicores, some of the higher level memories are shared between
  multiple cores. Moreover enforcing data locality has a major
  drawback in the fact that it seriously limits the scheduling of
  tasks since each core can only be assigned tasks that operate on
  data that resides on the memory associated with it. Preliminary
  results show that enforcing data locality and CPU affinity provides
  a slight speedup on the 8-way Dual Opteron system which is a NUMA
  architecture. These techniques require further investigation.

\item[Implement the same algorithms in distributed memory
  systems.] The fact that the block algorithms for QR and LU
  factorizations only require
  loose synchronization between tasks make them also good candidates
  for the implementation on distributed memory systems based on MPI
  communications.

\item[Implement the same algorithms on the STI Cell architecture.] In
  the STI Cell processor, no caches are present but a small,
  explicitly managed memory is associated to each core. Due to the
  small size of these local memories (only 256 KB), the LAPACK
  algorithms for LU and QR factorizations cannot be efficiently
  implemented. The block algorithms for LU and QR factorizations
  represent ideal candidates for the STI Cell architecture since they
  can be parallelized with a very fine granularity.

\item[Explore the usage parallel programming environments.]  The task of
  implementing Linear Algebra operations with
  dynamic scheduling of tasks on multicore architectures can be
  considerably simplified by the use of graph driven parallel
  programming environments. One such environment is SMP Superscalar\cite{smpss}
  developed at the Barcelona Supercomputing Center. SMP
  Superscalar addresses the automatic exploitation of the functional
  parallelism of a sequential program in multicore and SMP
  environments. The focus in on the portability, simplicity and
  flexibility of the programming model. Based on a simple annotation
  of the source code, a source to source compiler generates the
  necessary code and a runtime library exploits the existing
  parallelism by building at runtime a task dependency graph. The
  runtime takes care of scheduling the tasks and handling the
  associated data. Besides, a temporal locality driven task scheduling
  can be implemented.

\end{description}

\section*{Acknowledgment}
We would like to thank John Gilbert from University of California at Santa Barbara
for sharing with us its 8-way dual Opteron system.

\bibliography{blk_qr}  
\bibliographystyle{plain}     

\end{document}